\documentclass[11pt]{amsart}

\usepackage{epsfig, color}

\usepackage{amsmath}
\usepackage{amsthm}
\usepackage{hyperref}

\newtheorem{prop}{Proposition}
\newtheorem{lemma}[prop]{Lemma}
\newtheorem{thm}[prop]{Theorem}

\theoremstyle{definition}

\newcommand{\del}{\partial}
\newcommand{\dt}{\frac{\partial}{\partial t}}
\newcommand{\brs}[1]{\left| #1 \right|}

\newcommand{\gD}{\Delta}
\newcommand{\gd}{\delta}

\newcommand{\gl}{\lambda}

\newcommand{\ga}{\alpha}

\renewcommand{\ge}{\epsilon}
\newcommand{\N}{\nabla}
\newcommand{\FF}{\mathcal F}

\newcommand{\til}[1]{\widetilde{#1}}

\newcommand{\ohat}[1]{\overset{\circ}{#1}}

\DeclareMathOperator{\Rm}{Rm}

\DeclareMathOperator{\Id}{Id}
\DeclareMathOperator{\divg}{div}
\DeclareMathOperator{\grad}{grad}
\DeclareMathOperator{\Vol}{Vol}

\DeclareMathOperator{\Ca}{Ca}

\begin{document}

\title[$L^2$ flow on surfaces]{The
gradient flow of the $L^2$ curvature energy on surfaces}
\author{Jeffrey Streets}
\address{Fine Hall\\
         Princeton University\\
         Princeton, NJ 08544}
\email{\href{mailto:jstreets@math.princeton.edu}{jstreets@math.princeton.edu}}

\thanks{The author was partly supported by the National Science Foundation
via
DMS-0703660}

\begin{abstract} We investigate the gradient flow of
the $L^2$ norm of the Riemannian curvature on surfaces.  We show long time
existence with arbitrary initial data, and exponential convergence of the volume
normalized flow to a constant scalar curvature metric when the initial energy is
below a constant determined by the Euler characteristic of the underlying
surface.
\end{abstract}

\date{July 10th, 2010}

\maketitle

\section{Introduction}

In this paper we study the gradient flow of the
$L^2$ norm of the curvature tensor on surfaces.  Let us
first introduce some notation.  Fix $(M^n, g)$ a Riemannian manifold, and let
$\Rm$ denote the Riemannian curvature
tensor, and $s$ the scalar curvature.  Furthermore, let
\begin{gather*}
\mathcal F(g) := \int_M \brs{\Rm_g}_g^2 dV_g.
\end{gather*}
In what follows we will often drop the explicit reference to $g$, as all objects
in sight will be referencing a given time-dependent metric.  A basic calculation
(\cite{Besse} Proposition 4.70) shows that
\begin{align} \label{eq:gradF}
\grad \FF =&\ \gd d r - \check{R} + \frac{1}{4} \brs{\Rm}^2 g.
\end{align}
where $d$ is the exterior derivative acting on the Ricci tensor treated
as a one-form with values in the cotangent bundle, and $\gd$ is the adjoint of
$d$.  Moreover, 
\begin{align*}
\check{R}_{ij} = R_{i p q r} R_j^{p q r}.
\end{align*}
For purposes of this paper, we say that a metric is \emph{critical} if
\begin{align*}
\grad \FF \equiv \gl g
\end{align*}
for some $\gl$.  In general gradient flows of Riemannian functionals are natural
tools for finding critical points.  Specifically, consider the negative gradient
flow of $\mathcal F$:
\begin{gather} \label{flow}
\begin{split}
\dt g =&\ - \grad \FF,\\
g(0) =&\ g_0.
\end{split}
\end{gather}
This is a nonlinear fourth order degenerate parabolic equation.  Since the
equation is fourth order maximum principle techniques are not available,
and the analysis largely relies on integral estimates.  In \cite{Streets} we
showed
short-time existence of the
initial value problem as well as derivative estimates and a long-time existence
obstruction.  Furthermore, in \cite{Streets2} we showed further regularity
results, and in \cite{Streets3} showed exponential convergence in four
dimensions to a round metric when the $L^2$ norm of the traceless curvature
tensor is sufficiently small.

Due to the scaling properties, the functional $\FF$ is not particularly
interesting
in dimensions $n \geq 5$, but is certainly interesting when $n \leq 4$.  From
here on let us fix $n = 2$.
In this case the energy $\mathcal F$ takes the simple form (after rescaling),
\begin{align*}
\mathcal F(g) = \int_M s_g^2 dV_g.
\end{align*}
By adding a constant depending only on the genus of $M$, setting
\begin{gather*}
\bar{s}_g = \frac{\int_M s_g dV_g}{\int_M dV_g}
\end{gather*}
we may express
\begin{align} \label{energy}
\mathcal F(g) + \frac{16 \pi^2 \chi(M)^2}{\Vol(M)} = \int_M \left(s_g -
\bar{s}_g \right)^2 dV_g =: \Ca(g),
\end{align}
where $\Ca(g)$ denotes the Calabi energy.  By applying the algebraic curvature
identity
\begin{align*}
R_{ijkl} = \frac{1}{2} s \left[ g_{il} g_{jk} - g_{ik} g_{jl} \right]
\end{align*}
which holds on surfaces one can compute from (\ref{eq:gradF}) that
\begin{align} \label{eq:gradFsurfaces}
\grad \FF =&\ - \gD s g + \N^2 s - \frac{1}{4} s^2 g.
\end{align}
Thus we may express equation (\ref{flow}) in the simpler form
\begin{align} \label{surfacesflow}
\dt g =&\ \gD s g - \N^2 s + \frac{1}{4} s^2 g.
\end{align}
We will also consider solutions to the corresponding volume-normalized equation
\begin{align} \label{vnflow}
\dt g =&\ \gD s g - \N^2 s + \frac{1}{4} s^2 g - \frac{1}{4} \frac{\int_M s^2
dV}{\int_M dV} g
\end{align}
Critical metrics on surfaces, i.e. metrics where the right hand side of
(\ref{vnflow}) vanishes, are also called \emph{extremal Hermitian metrics}, and
have been the object of intense study (see for instance \cite{Chenextremal}). 
Indeed,
Chen's convergence/compactness criterion for metrics with bounded volume and
Calabi energy \cite{Chen} will be essential to our proofs.

\begin{thm} \label{surfaces} Let $(M^2, g_0)$ be a compact surface.  The
solutions to (\ref{surfacesflow}) and (\ref{vnflow}) with initial condition
$g_0$ exist
for all time.
\end{thm}

\noindent Incidentally, equation (\ref{vnflow}) is the gradient flow of the
scale-invariant version of $\FF$, i.e.
\begin{align*}
 \mathcal E(g) =&\ \FF(g) \Vol(g).
\end{align*}
Observe that, as a consequence of the Gauss-Bonnet theorem and H\"older's
inequality, we have
\begin{align*}
\left[4 \pi \chi(M) \right]^2 =&\ \left(\int_M s dV \right)^2 \leq \mathcal
E(g),
\end{align*}
with equality achieved only in the case of constant scalar curvature metrics. 
If the initial energy is sufficiently close to this minimum value, we are able
to conclude exponential convergence of (\ref{vnflow}) at infinity. 
Specifically,

\begin{thm} \label{convergence} Let $(M^2, g_0)$ be a compact surface.  If
\begin{align} \label{lowenergy}
\mathcal E(g_0) < 16 \pi^2 \left[ \brs{\chi(M)} + 1 \right]^2.
\end{align}
the solution to (\ref{vnflow}) converges exponentially to
a constant scalar curvature metric.
\end{thm}

\section{Proof of Theorem \ref{surfaces}}
We will start the proof of Theorem \ref{surfaces} with a reduction.  Let
$(M^2, g(t))$ be a solution to (\ref{surfacesflow}).  The Hessian term can be
expressed
as a Lie derivative, specifically $\N^2 s =
\frac{1}{2} \mathcal L_{\N s} g$.  Therefore by a well-known procedure we can
modify a solution to (\ref{flow}) by a family of diffeomorphisms
to remove this term.  We will make this procedure more precise below.  For now
consider the equation
\begin{align} \label{conformalflow}
\dt g =&\ \gD s g + \frac{1}{4} s^2 g.
\end{align}
Equation (\ref{conformalflow}) has the nice advantage of being conformal.  It
follows from the diffeomorphism invariance of $\mathcal F$ that it is
nonincreasing for solutions to (\ref{conformalflow}).  However, we need the
precise form the evolution equation.

\begin{lemma} \label{surfaceslemma1} Let $(M^2, g(t))$ be a solution to
(\ref{conformalflow}).  Then
\begin{align*}
 \dt \mathcal F(g(t)) =&\ - \int_M \left[\frac{1}{2} \left(\gD s + \frac{1}{2}
s^2\right)^2 + \brs{\ohat{\N^2 s}}^2 \right]dV
\end{align*}
where $\ohat{\N^2 s}$ is the component of $\N^2 s$ orthogonal to $g$.
\begin{proof} Using (\ref{eq:gradFsurfaces}) we directly compute that a solution
to (\ref{conformalflow}) satisfies
\begin{align*}
\dt \mathcal F(g(t)) =&\ \int_M \left< - \gD s g + \N^2 s - \frac{1}{4} s^2 g,
\gD s g + \frac{1}{4} s^2 g \right> dV
\end{align*}
However, $\grad \FF$ is divergence-free, therefore
\begin{align*}
\int_M \left< - \gD s + \N^2 s - \frac{1}{4} s^2 g, \N^2 s \right> dV =&\ \int_M
\left< \divg \grad \FF, \N s \right> = 0.
\end{align*}
Hence
\begin{align*}
\dt \mathcal F(g(t)) =&\ - \int_M \brs{ \gD s g - \N^2 s + \frac{1}{4} s^2 g}^2
dV.
\end{align*}
Now let us decompose pointwise
\begin{align*}
\gD s g - \N^2 s g + \frac{1}{4} s^2 g =&\ \gD s g - \frac{1}{2} \gD s g -
\ohat{\N^2 s} + \frac{1}{4} s^2 g\\
=&\ \frac{1}{2} \left( \gD s + \frac{1}{2} s^2 \right) g - \ohat{\N^2 s}
\end{align*}
where $\ohat{\N^2 s}$ refers to the component of $\N^2 s$ orthogonal to $g$. 
Therefore one has, pointwise,
\begin{align*}
\brs{ \gD s g - \N^2 s + \frac{1}{4} s^2 g}^2 =&\ \brs{ \frac{1}{2} \left( \gD s
+ \frac{1}{2} s^2 \right) g - \ohat{\N^2 s}}^2\\
=&\ \brs{ \frac{1}{2} \left( \gD s + \frac{1}{2} s^2 \right) g}^2 +
\brs{\ohat{\N^2 s}}^2\\
=&\ \frac{1}{2} \left(\gD s + \frac{1}{2} s^2 \right)^2 + \brs{\ohat{\N^2 s}}^2.
\end{align*}
The result follows.
\end{proof}
\end{lemma}

We recall Chen's bubbling criterion for conformal metrics on surfaces.  In the
statement below, $B_R(x)$ refers to the ball of radius
$R$ around the point $x$ in a \emph{background} metric $g_0$.

\begin{thm} \label{chen} (\cite{Chen} Theorem 1, \cite{Struwe} Theorem 3.2). Let
$g_n = e^{2 u_n} g_0$ be a sequence of smooth conformal metrics on $M$ with unit
volume and bounded Calabi energy.  Then either the sequence $\{u_n\}$ is bounded
in $H^2(M, g_0)$ or there exist points $x_1, \dots x_l \in M$ and a subsequence
$\{u_n\}$  such that for any $R > 0$ and any $l \in \{1, \dots, L\}$ there holds
\begin{align*}
\liminf_{n \to \infty} \int_{B_R(x_l)} \brs{s_n} dV_{g_n} \geq 4 \pi.
\end{align*}
Moreover, there holds
\begin{align*}
4 \pi L \leq \limsup_{n \to \infty} \mathcal F(g_n)^{\frac{1}{2}} < \infty,
\end{align*}
and either $u_n \to - \infty$ as $n \to \infty$ locally uniformly on $M
\backslash \{x_1, \dots, x_L\}$ or $\{u_n\}$ is locally bounded in $H^2(M, g_0)$
away from $x_1, \dots x_L.$
\end{thm}

\noindent We can now give the proof of Theorem \ref{surfaces}:
\begin{proof} Given $(M^2, g)$, let $g(t)$ be the solution to
(\ref{conformalflow}) with initial condition $g$.   By Lemma
\ref{surfaceslemma1} $\Ca(g(t)) \leq \Ca(g)$.  Suppose $T < \infty$ is the
maximal existence time of $g(t)$. 
According to Theorem \ref{chen} we either have that $u_t$ is uniformly bounded
in $H^2(M, g_0)$ as $t \to T$ or there is at least one point $x \in M$ such that
for all $R > 0$,
\begin{align} \label{surfacesloc10}
 \lim_{t \to T} \int_{B_R(x)} \brs{s_t} \geq 4 \pi.
\end{align}
Suppose the second alternative holds, and let $x$ be the point where
(\ref{surfacesloc10}) holds.  Fix some $\ge > 0$.  Since $\mathcal F$ is a
continuous,
nonincreasing function of $t$ which is bounded below, we may choose $t_0 > 0$,
$T - 1 < t_0 < T$ such that
\begin{align*}
\mathcal F(t_0) - \lim_{t \to T} \mathcal F(t) \leq \ge.
\end{align*}
Furthermore, choose $R > 0$ small such that
\begin{align*}
\int_{B_R(x)} dV_{g({t_0})} \leq \ge.
\end{align*}
Let
\begin{align*}
 f(t) = \int_{B_R(x)} dV_{g(t)}.
\end{align*}
Recall that the integral is over the ball of radius $R$ with respect to $g_0$. 
We estimate
\begin{align*}
 \frac{d}{dt} f =&\ \int_{B_R(x)} \left( \gD s + \frac{1}{4} s^2 \right) dV_g\\
=&\ \int_{B_R(x)} \left( \gD s + \frac{1}{2} s^2  - \frac{1}{4} s^2 \right)
dV_g\\
\leq&\ \int_{B_R(x)} \left(1 + \left(\gD s + \frac{1}{2} s^2 \right)^2 \right)
dV\\
\leq&\ f(t) - 2 \frac{d}{dt} \mathcal F(g(t)).
\end{align*}
The last line follows using Lemma \ref{surfaceslemma1}.  We may integrate the
limiting ODE of this differential inequality to conclude that for $t \geq t_0$,
\begin{align*}
f(t) \leq \left( \int_{t_0}^t \left(-2\frac{d}{dp} \mathcal F(g(p)) \right)
e^{-(p - t_0)} dp + \ge \right) e^{t - t_0}.
\end{align*}
Since $ -\frac{d}{dp} \mathcal F(g(p)) \geq 0$ we can estimate
\begin{align*}
\int_{t_0}^t \left( -2 \frac{d}{dp} \mathcal F(g(p)) \right) e^{-(p - t_0)} dp
\leq&\ \int_{t_0}^t \left( -2 \frac{d}{dp} \mathcal F(g(p)) \right) dp\\
=&\ 2 \left(\mathcal F(g(t_0)) - \mathcal F(g(t)) \right)\\
\leq&\ 2 \left( \mathcal F(g(t_0)) - \lim_{t \to T} \mathcal F(g(t)) \right)\\
\leq&\ 2 \ge.
\end{align*}
It follows that
\begin{align*}
f(t) \leq 3 \ge e^{t - t_0}.
\end{align*}
We conclude that
\begin{align*}
 \lim_{t \to T} \int_{B_R(x)} \brs{s} dV_g \leq&\ \lim_{t \to T} \left(
\int_{B_R(x)} s^2 dV_g \right)^{\frac{1}{2}} f(t)^{\frac{1}{2}}\\
\leq&\ \left( 3 \ge e^{T - t_0} \mathcal F(g(0))\right)^{\frac{1}{2}}.
\end{align*}
Since $t_0 > T - 1$, we conclude that if $\ge \leq\frac{1}{3 e \mathcal
F(g(0))}$, then
\begin{align*}
\lim_{t \to T} \int_{B_R(x)} \brs{s} dV_g < 4 \pi.
\end{align*}
for $\ge$ chosen small with respect to $\mathcal F(g(0))$.  This contradicts
(\ref{surfacesloc10}), and therefore the first alternative of Theorem \ref{chen}
must hold.  In particular, 
$u_t$ is uniformly bounded in $H^2$, and hence in $C^0$ by the Sobolev
embdedding.  We next claim that the curvature is bounded up to time $T$. 
Suppose
\begin{align*}
 \limsup_{t \to T} \brs{s} = \infty
\end{align*}
and choose a sequence $\{x_i, t_i\}$ realizing this limit, and let $\gl_i =
\brs{s}_{g(t_i)}(x_i)$.  Consider the sequence of solutions to
(\ref{conformalflow}),
\begin{align*}
g_i(t) := \gl_i g \left( t_i + \frac{t}{\gl_i^2} \right).
\end{align*}
By construction $\sup \brs{s}_{g_i}(0) = 1$.  The sequence of metrics $\{g_i\}$
have bounded Sobolev constant (which is scale invariant) and bounded curvature. 
It follows from (\cite{Streets} Theorem 7.1) that there exists a subsequence
converging to a limiting solution $g_{\infty}$ of (\ref{conformalflow}). 
Strictly speaking (\cite{Streets} Theorem 7.1) was only proved for solutions
to (\ref{flow}), but the estimates clearly apply since they are diffeomorphism
invariant.  By construction $\sup \brs{s}_{g_{\infty}} = 1$. 
However, $\FF(g(t_i)) \leq \FF(g(0))$, and moreover $\FF(\gl g) = \frac{1}{\gl}
\FF(g)$.  Thus $\FF(g_{\infty}) = 0$, a contradiction.  It follows that the
curvature is bounded up to time $T$, and hence by (\cite{Streets} Theorem
6.2), the solution exists smoothly past time $T$.

We use a diffeomorphism change to obtain the long-time existence of solutions of
(\ref{flow}).  Let $g(t)$ denote the solution to (\ref{conformalflow}) with
initial condition $g$ discussed above.  Following the prior discussion, we may
define the one-parameter family of diffeomorphisms
\begin{align*}
\dt \phi_t =&\ \frac{1}{2} \left( \N s \right)^{\sharp}\\
\phi_0 =&\ \Id
\end{align*}
where here $(\N s)^{\sharp}$ denotes the vector field dual to $\N s$ with
respect to
the metric $g(t)$.  By the estimates for $g(t)$ which we have shown, $\phi_t$ is
smooth on any finite time interval $[0, T]$.  Consider the one-parameter family
of metrics $\til{g}(t) = \phi_t^* g(t)$.  Note that by construction $\til{g}(0)
= g(0) = g$.  We can directly compute
\begin{align*}
\dt \til{g}(t) =&\ \phi_t^* \left( \dt g \right) + \frac{\del}{\del s} \phi_{t +
s}^* g(t)\\
=&\ \left( - \grad \FF - \N^2 s \right)(\phi^*_t g(t)) - \mathcal
L_{(\phi_t^{-1})_* \N s(t)} \phi_t^* g(t))\\
=&\ - \grad \FF(\til{g}(t)) - \N^2 s(\phi_t^* g(t)) + \N^2 s(\phi_t^* g(t))\\
=&\ - \grad \FF(\til{g}(t)).
\end{align*}
Thus $\phi_t^*(g(t))$ is the solution to (\ref{flow}) with initial condition
$g$, which by the estimates already shown above exists for all time.
\end{proof}

\section{Proof of Theorem \ref{convergence}}

\begin{proof}
Let $g(t)$ be a solution to
the conformal volume-normalized flow, i.e.
\begin{align} \label{vnflow2}
\dt g =&\ \gD s g + \frac{1}{4} s^2 g - \frac{1}{4} \frac{\int_M s^2 dV}{\int_M
dV} g.
\end{align}
This flow preserves the initial volume, which we will assume without loss of
generality to be $1$.  This flow differs from (\ref{conformalflow}) by a
rescaling in space and a reparameterization of time, and these scales are
determined by the volume of the solution to (\ref{conformalflow}), which is
easily seen to grow at worst linearly.  In particular, it is easy to see that
the corresponding solution to (\ref{vnflow2}) exists for all time.  Furthermore,
one can show using an argument akin to Lemma \ref{surfaceslemma1} that
\begin{gather} \label{vngradient}
\begin{split}
\dt \mathcal F =&\ - \int_M \brs{\gD s g - \N^2 s + \frac{1}{4} s^2 g -
\frac{1}{4} \left(\int_M s^2 dV \right)g}^2\\
=&\ - \int_M \left[ \frac{1}{2} \left( \gD s + \frac{1}{2} s^2 - \frac{1}{2}
\left( \int_M s^2 dV \right)\right)^2 + \brs{\ohat{\N^2 s}}^2 \right] dV.
\end{split}
\end{gather}

First we want to show subsequential convergence to a metric of constant
curvature.  Consider now some sequence of times $\{t_j\} \to \infty$.  The
metrics $g(t_i)$ have bounded Calabi energy and unit volume.  By Theorem
\ref{chen} , we either have bubbling or smooth convergence.  In the case of
smooth convergence, the limiting metric is necessarily extremal Hermitian by
(\ref{vngradient}), and hence constant curvature.  This was first
observed by Calabi \cite{Calabi}, and can be seen by pairing the critical equation 
equation against $\gD s$ and integrating by parts.  So,
let us assume that the sequence bubbles in the sense of Theorem \ref{chen}.  At
this point we break into two cases.

First assume $\chi(M) \leq 0$.  By Theorem \ref{chen}, either one gets a
well-defined limit metric $g_{\infty}$ on $\til{M} := M \backslash \{\mbox{
bubble points}\}$ or the sequence experiences collapse and the volume of
$\til{M}$ is zero in the limit.  First we rule out the collapsing possibility. 
In the terminology of \cite{Chen}, if this occurs then $\til{M}$ is a ``ghost
vertex'' in the tree decomposition of the limit of $\{g(t_i)\}$.  Every vertex
which is not a ghost vertex (and one must exist since the area is fixed)
corresponds to a finite-area metric which is critical by (\ref{vngradient}). 
This is necessarily the rotationally symmetric ``teardrop metric'' of
Chen (\cite{Chenextremal} section 8).  One can show that this metric still has
strictly positive curvature integral, in particular, the contribution to the
Euler characteristic of the limiting tree is positive.  It follows that the
Euler characteristic of $\til{M}$ must be strictly negative.  We conclude that
\begin{align*}
C \leq \lim_{i \to \infty} \left( \int_{\til{M}} K_{g_i} dV_i \right)^2 \leq
\left(\int_{\til{M}} dV_i \right)
\left( \int_{\til{M}} K^2 dV_i \right). 
\end{align*}
But the area approaches zero and the Calabi energy remains finite, so this is a
contradiction.  Therefore the collapse possiblity does not occur, and hence we
conclude weak convergence to a critical bubbled metric.  The schematic picture
of the limit is in Figure
\ref{fig:figure1}.

\begin{figure}[ht]
\begin{center} \resizebox{275pt}{!}{\input{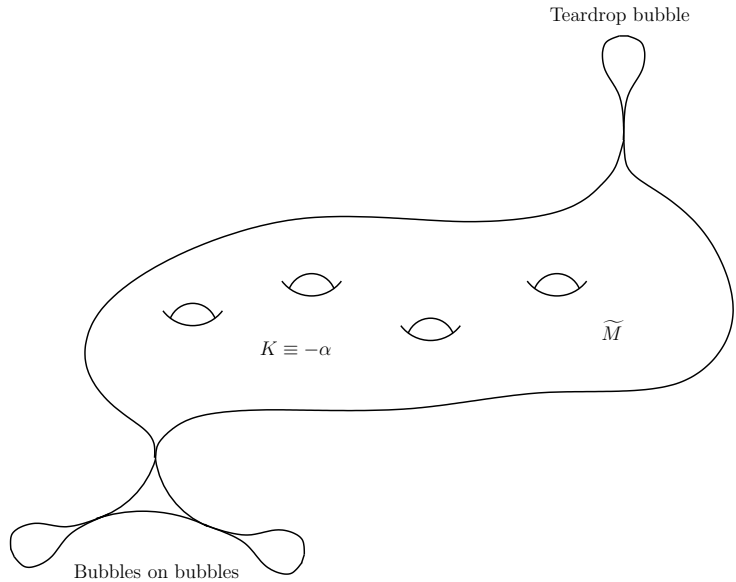}}
\end{center}
\caption{A conformal bubbled metric} \label{fig:figure1}
\end{figure}

It remains to compute the energy of this metric to derive a contradiction. 
Since there exists at least one bubble point, the Euler characteristic of
$\til{M}$ is strictly negative.  Specifically,
\begin{align*}
\chi(\til{M}) \leq \chi(M) - 1.
\end{align*}
Since $g_{\infty}$ is an extremal metric on a surface of negative Euler
characteristic, it follows from (\cite{Chenextremal} Theorem A) that in fact
$s_{g_{\infty}} \equiv - \ga$.  We conclude by the Gauss-Bonnet theorem that
\begin{align*}
4 \pi \brs{\chi(\til{M})} = \brs{\int_{M_{\infty}} s dV} = \ga \Vol(\til{M}).
\end{align*}
Therefore since $\Vol(\til{M}) \leq 1$, we directly estimate
\begin{align*}
\mathcal F(g_{\infty}) =&\ \int_{M_{\infty}} s^2 dV\\
=&\ \ga^2 \Vol(M_{\infty})\\
=&\ \frac{16 \pi^2 \brs{\chi(\til{M})}^2} {\Vol(M_{\infty})}\\
\geq&\ 16 \pi^2 \left[ \brs{\chi(M)} + 1 \right]^2.
\end{align*}
Since the volume was fixed to be $1$ along the flow, we have by hypothesis that
the initial energy is below $16 \pi^2 \left[ \brs{\chi(M)} + 1 \right]^2$. 
Certainly then the limit should have energy below this value as well,
contradicting the above.  Thus this bubbling possiblity does not occur and the
flow converges to a constant scalar curvature metric, as claimed. 

Now we address the case of the sphere.  First let us rule out the collapsing
case.  By following the proof in the negative Euler characteristic case we see
that the only case which is not immediately ruled out is when the Euler
characteristic of $\til{M}$ (using the notation above) is zero, which can
theoretically happen in the case of two bubbling points.  The simplest example
of this behavior is shown in Figure \ref{fig:figure2}.

\begin{figure}[ht]
\begin{center} \resizebox{275pt}{!}{\input{Collapsedsphere.pstex_t}}
\end{center}
\caption{Collapsed metric on the sphere} \label{fig:figure2}
\end{figure}

In this case, there must be at least two ``deepest bubbles'' $B_1, B_2$ which by
the arguments above must correspond to teardrop metrics.   One can explicitly
compute the energy contained in a teardrop bubble in terms of its area.  In
particular (see \cite{Chenextremal} section 8), one has
\begin{align*}
\int_{B_i} s^2 = \frac{16 \pi^2}{\Vol(B_i)}.
\end{align*}
Since $\Vol(B_1) + \Vol(B_2) \leq 1$, one can use the inequality $\frac{1}{x} +
\frac{1}{1-x} \geq 4$ for $0 < x < 1$ to conclude that
\begin{align} \label{sphereenergy}
\lim_{i \to \infty} \mathcal F(g(t_i)) \geq 64 \pi^2,
\end{align}
contradicting our low-energy hypothesis.  Therefore we must get a noncollapsed
limit at infinity, and it remains to rule out bubbling in this case.  We know
that in this noncollapsed case that the limiting metric must be a critical
metric with cusp ends.  In particular, there must be at least three bubble
points and the limiting metric has constant negative curvature.  What is more,
there are three deepest bubbles which must correspond to teardrop metrics, and
hence precisely as above one obtains (\ref{sphereenergy}).  In fact one gets an
even better bound, but certainly at least we have proved that in the general
case (\ref{sphereenergy}) holds, contradicting our initial energy hypothesis.

To finish the proof we need to exhibit exponential convergence, which we only
sketch here since these details are the same as the case of Calabi flow which
appear in (\cite{Chencalabi} Proposition 6.1).  The main estimate to show is
exponential decay of the Calabi energy.  By the above arguments, any sequence of
times $\{t_i \}$ must be bounded in $H^2$, therefore there is a uniform bound
for $u$ in $H^2$ which holds for all times.  By the Sobolev embedding we
conclude a uniform $C^0$ bound for $u$, and hence a uniform bound for the
Sobolev constants of the time-varying metrics.  Moreover, the Calabi energy
converges subsequentially to zero.  At any large time where $\Ca(g(t)) \leq
\ge$, using the Sobolev constant estimate of the time-varying metrics, one can
obtain the estimate
\begin{align*}
\dt \Ca(g) \leq&\ -(1 - C \ge) \int_M \brs{\gD (s - \bar{s})}^2 dV_g +
\frac{\bar{s}}{8} \int_M \brs{\N s}^2 dV_g.
\end{align*}
In the case of nonpositive Euler characteristic one immediately concludes
exponential convergence.  In the case of the sphere one has to exploit the
Kazdan-Warner identity to show that $s - \bar{s}$ is nearly orthogonal to the
first eigenspace of $\gD$, after which exponential convergence follows.  With
this exponential decay in place we can get exponential convergence in stronger
norms.  First, since we know the Sobolev constant is uniformly bounded for all
time, by repeating the blowup argument of Theorem \ref{surfaces} we obtain a
uniform curvature bound.  Now we can apply (\cite{Streets} Theorem 5.4) to
conclude that all $H_k^2$ norms of curvature decay exponentially as well.  With
the exponential convergence in hand, standard calculations show that once we
pull back the conformal flow to give a solution of (\ref{flow}), we still have
exponential convergence, and the result follows.
\end{proof}

\section{Further Questions}

Naturally, one expects that the flow should converge to a constant scalar
curvature metric with arbitrary initial data.  The reason one can obtain this
for Calabi flow is that one has monotone quantities which are constant only at
constant scalar curvature metrics.  Indeed, the existence of such a quantity for
(\ref{flow}) would immediately imply convergence with arbitrary initial data,
via arguments that appear in (\cite{Chencalabi}).

\bibliographystyle{hamsplain}

\end{document}